\documentclass[10pt,runningheads,a4paper]{llncs}
\usepackage{graphicx}
\usepackage{graphicx}
\usepackage{amsmath}
\usepackage{amscd}

\newtheorem{open}{Open Problem}
\begin{document}

\title{On cycle-supermagic labelings of the disconnected graphs}
\author{Syed Tahir Raza Rizvi, Kashif Ali
 }

\institute{ Graphs and Combinatorics Research Group,\\
Department of Mathematical Sciences,\\
COMSATS Institute of Information Technology, \\
Lahore, Pakistan \\ \{ strrizvi, akashifali@gmail.com\}}
\authorrunning{ Rizvi et al.}
\titlerunning{On cycle-supermagic labelings of the disconnected graphs}
\maketitle

\begin{abstract}
A graph $G(V,E)$ has an $H$-covering if every edge in $E$ belongs to
a subgraph of $G$ isomorphic to $H$. Suppose $G$ admits an
$H$-covering. An $H$-magic labeling is a total labeling $\lambda$
from $V(G) \cup E(G)$ onto the integers $\{1, 2, \cdots ,|V(G) \cup
E(G)|\}$ with the property that, for every subgraph $A$ of $G$
isomorphic to $H$ there is a positive integer $c$ such that $\sum
A=\sum\limits_{v\in V(A)}\lambda(v) + \sum\limits_{e\in
E(A)}\lambda(e)= c .$ A graph that admits such a labeling is called
{\em $H$-magic}. In addition, if $\{\lambda(v)\}_{v\epsilon
V}=\{1,2, \cdots ,|V|\}$, then the graph is called $H$-supermagic.
In this paper we formulate cycle-supermagic labelings for the
disjoint union of isomorphic copies of different families of graphs.
We also prove that disjoint union of non isomorphic copies of fans
and ladders are cycle-supermagic.

\end{abstract}

\begin{flushleft}
\end{flushleft}

\section{Introduction}
Let $G=(V,E)$ be a finite, simple, planar, connected and undirected
graph, where $V$ and $E$ are its  vertex-set and edge-set
respectively.  A {\em labeling} (or {\em valuation}) of a graph is a
map that carries graph elements to numbers (usually positive or
non-negative integers). Let $H$ and $G=(V,E)$ be
finite simple graphs with the property that every edge of $G$
belongs to at least one subgraph isomorphic to $H$. A bijection $\lambda:
{V\cup E }\rightarrow\{1,...,|V|+|E|\}$ is an $H$-magic labeling of
$G$ if there exist a positive integer $c$ (called the magic constant),
such that for any subgraph $H'(V', E')$ of $G$ isomorphic to $H$,
the sum
$\sum\limits_{v\in V'}\lambda(v) + \sum\limits_{e\in E'}\lambda(e)$
is equal to magic constant $c$. This sum is also known as weight of $H$. A graph is $H$-magic if it admits
an $H$-magic labeling. In addition, if $H$-magic labeling $\lambda$
has the property that $\{\lambda(v)\}_{v\epsilon V}=\{1, \cdots ,|V|\}$,
then the graph is $H$-supermagic. In our terminology, super edge-magic total labeling
is a $K_2$-supermagic labeling. The notion of $H$-magic graph was introduced by
A.Gutierrez and A. Llado [11],  as an extension of the magic
valuation given by Rosa \cite{KR06}, which corresponds to the case
$H \cong K_2.$ \\

\noindent There are many results for $H$-(super)magic graphs. Llado and Moragas \cite{KR066}
 studied cycle-magic labeling of wheels, prisms, books, wind
mill graphs and subdivided wheels by using the technique of partitioning sets of integers.  Ngurah et al. \cite{FIM031} constructed
cycle-supermagic labelings for some families of graphs namely fans,
ladders, and books etc. For any connected graph $G$, Maryati et al. \cite{MSB08},
 proved that disjoint union of $k$ isomorphic copies of $G$ is a $G$-supermagic
graph if and only if $|V(G)|+|E(G)|$ is even or $k$ is odd. For a disconnected graph $G,$ having at least 2 vertices in its every components, the result was proved  again by
Maryati et al. \cite{FIM03}. For a positive integer $k$, they proved that (i) If $kG$ is $G$-magic, then $|V(G)|+|E(G)|$ is even or $k$ is odd and (ii) If $|V(G)|+|E(G)|$ is even or $k$ is odd, then $kG$ is $G$-supermagic.
For further detail, see the recent survey by Gallian \cite{G05}.\\

\noindent In this paper, we study the problem that if a connected graph $G$  is cycle-supermagic, then either disjoint
 union of $m$ isomorphic copies of graph $G$ denoted by $mG,$ is cycle-supermagic or not? We study cycle-supermagic labelings for the disjoint
union of isomorphic copies of fans, ladders, triangular ladders,
wheel graphs, book graph and generalized antiprism graph. Also, we
study cycle-supermagic labelings of disjoint union of non-isomorphic
copies of ladders and fans.

\section{Main Results}

\noindent Before giving our main results, let us consider the following lemma
found in \cite{KR06} that gives a necessary and sufficient
condition for a graph to be super edge-magic.
\begin{lemma} \label{L01}
A graph $G$ with $v$ vertices and $e$ edges is super edge-magic
 if and only if there exists a bijective function $\lambda :V(G)
\rightarrow \{1, 2, \cdots, v\}$ such that the set $S = \{\lambda(x) +
\lambda(y) | xy \in E(G)\}$ consists of $e$ consecutive integers. In such
a case, $\lambda$ extends to a super edge-magic labeling of $G$ with
magic constant $c = v + e + s$, where $s = min(S).$
\end{lemma}
\vspace{.5cm}
\noindent Moreover, while proving our some results,
we will use following result to formulate
the $(a,2)$-edge-antimagic vertex labeling for
 a subgraph isomorphic to $mP_n$.

\begin{theorem}\cite{Dafik09} The graph $mP_n$ has a $(m+2,2)$-edge-antimagic vertex labeling, for every $m\geq2$ and $n\geq2$.

\end{theorem}

\subsection{Cycle-supermagic labelings of the disjoint union of isomorphic graphs}
In the following theorem, we consider $C_3$-supermagic
labeling for disjoint union of isomorphic copies of fans. For $n\geq3,$ the {\em fan} $F_n\cong K_1+P_n$ is a graph with $V(F_n)
=\{c\}\cup\{v_i:1\leq i \leq n\}, $ $E(F_n) =\{cv_i:1\leq i \leq
n\}\cup  \{v_iv_{i+1}:1\leq i \leq n-1\}. $

\begin{theorem}
For any positive integer $m\geq2,$ and $n\geq 3$, the graph $G\cong mF_n$
is $C_3$-supermagic.
\end{theorem}
\noindent{\bf Proof.}\\
Let $v=|V(G)|$ and $e=|E(G)|$. Then $v=m(n+1)$, $e=m(2n-1)$. We
denote the vertex and edge sets of $G$ as follows:\\
\noindent $V(G)=\{c_j: 1\leq j\leq m \}\cup \{v^{j}_{i}: 1\leq i\leq n , 1\leq j\leq m \},$\\
\noindent $E(G)=\{v^{j}_{i}v^{j}_{i+1}: 1\leq i\leq n-1 , 1\leq j\leq m \}
\cup \{c_jv^{j}_{i}:1\leq i\leq n , 1\leq j\leq m\}. $\\
\noindent For $1\leq i \leq n-1, 1\leq j \leq m,$ let
$C^{(i,j)}_{3}$ be subcycle of graph $G$ with $V(C^{(i,j)}_{3})=\{c_j,v_i^j,v_{i+1}^j\}$ and $E(C^{(i,j)}_{3})=\{c_jv_{i}^j, v_i^jv_{i+1} ^ j,v_{i+1} ^ jc_j\}.$\\

\noindent
 We will define a total labeling $\lambda:V(G)\cup E(G)\rightarrow
\{1,2,....,3mn\}$ as follows:\\

\noindent Let $1\leq
i\leq n, 1\leq
j\leq m.$\\

\noindent {  Step 1.} At this step, we will only label the vertices of graph $G$ as follows:
$$\lambda(c_j)=m(n+1)-j+1$$
\noindent For $i$-odd

$$\lambda (v_i^j)=j+m(\frac{i-1}{2})$$

\noindent For $i$-even
$$
\lambda (v_i^j)= \frac{[2n+2i+3+(-1)^{n+1}]m+4j}{4}$$
After labeling the vertices of graph $G$, we can see that the set of all weights of subcycle
$C^{(i,j)}_{3}$ under labeling $\lambda$ consists of consecutive integers $\frac{3mn}{2}+m+2,\frac{3mn}{2}+m+3,...,\frac{5mn}{2}+1,$ for even $n$ and
$\frac{3m(n+1)}{2}+2,\frac{3m(n+1)}{2}+3,...,\frac{m(5n+1)}{2}+1,$ for odd $n$.\\

\noindent { Step 2.}
 Now, we will label the edges $c_jv_{i}^j$  as follow:
  $$\lambda(c_jv_{i}^j)=3mn-m(i-1)-j+1$$
  \noindent After labeling the vertices and edges $c_jv_{i}^j$ of graph $G$, we can see that the set of all weights of subcycle
$C^{(i,j)}_{3}$ consists of consecutive integers $\frac{m(26n+5+(-1)^{n+1})}{4}+3,\frac{m(26n+5+(-1)^{n+1})}{4}+4,...,\frac{m(26n+5+(-1)^{n+1})}{4}+m(n-1)+2,$ for any $n$.\\

\noindent { Step 3.} Now we are left with the labeling of edges $v_i^jv_{i+1} ^ j: 1\leq
i\leq n-1,$ with consecutive
integers $mn+m+1, mn+m+2,...,2mn.$
Here, we can easily label these edges by using Lemma 1 to have $C_3$-supermagic labeling of graph $G$ as follows:

$$\lambda(v_i^jv_{i+1} ^ j)=m(n+i)+j.$$

\noindent After combining all three steps, we can see that for every subcycle $C^{(i,j)}_{3}$ of $G$, the sum $\lambda (c_j)+ \lambda
(v_i^j)+\lambda
(v_{i+1}^j)+\lambda(c_jv_{i}^j)+\lambda(c_jv_{i+1}^j)+\lambda(v_i^jv_{i+1}
^ j)$ is
$\frac{m}{4}[34n+5+(-1)^{n+1}]+3.$ Hence $mF_n$ is $C_3$-supermagic.\qed\vspace{.5cm}
\noindent In the next two theorems, we prove that disjoint union of isomorphic copies of ladder graphs, namely ladder graphs and triangular
ladders   admit
cycle-supermagic labelings.\\
\noindent  Let $L_n\cong P_n\times P_2,$ $n\geq2,$ be a ladder graph with
$V(L_n)=\{u_i, v_i:\hphantom{a}1\leq i \leq n\},$
$E(L_n)=\{u_iv_i:1\leq i \leq
n\}\cup\{u_iu_{i+1}: 1\leq i \leq n-1\}\cup \{ v_i v_{i+1}:1\leq i \leq n-1\}$.\\

\begin{theorem}
For every $m\geq2$, $n\geq 2,$ the graph  $G\cong mL_n$ is
$C_4$-supermagic.
\end{theorem}
\noindent{\bf Proof.}\\
Let $v=|V(G)|=2mn,$ and $e=|E(G)|=3mn-2m$. We denote the vertex and edge sets of $G$ as follows:\\
\noindent
$V(G)=\{u^{j}_{i}: 1\leq i\leq n
, 1\leq j\leq m \}\cup \{ v^{j}_{i}: 1\leq i\leq n , 1\leq j\leq m
\};$\\ \noindent
 $E(G)=\{u^{j}_{i} v^{j}_{i}: 1\leq i\leq
n , 1\leq j\leq m \} \cup  \{u^{j}_{i}u^{j}_{i+1}:1\leq i\leq n-1 ,
1\leq j\leq m \}\cup \{ v^{j}_{i}v^{j}_{i+1}: 1\leq i\leq n-1 ,
1\leq j\leq m \}.$\\
\\

\noindent
 Define a total labeling $\lambda:V(G)\cup E(G)\rightarrow
\{1,2,...,5mn-2m\}$ as follows:\\

\noindent { Step 1.} For $1\leq i\leq n , 1\leq j\leq m$
$$\lambda (u_i^j)=j+m(i-1)$$
$$\lambda (v_i^j)=v-m(i-1)-j+1$$
After labeling the vertices of graph $G$, we can see that weight of every subcycle
$C^{(i,j)}_{4}$ under labeling $\lambda$ is $2v+2.$\\

\noindent {  Step 2.} Now, we will use Theorem 1 to label the edges $u_i^jv_{i} ^ j,$
 by assuming that
these edges are vertices of subgraph isomorphic to $mP_n$ and labeling is as follows:
$$\lambda(u_i^jv_{i} ^ j)=v+m(i-1)+j$$
At this stage, set of all weights of $C^{(i,j)}_{4}$ under labeling $\lambda$ forms
an arithmetic sequence with initial term $4v+m+2$ and the difference is 2.\\

\noindent { Step 3.}  For $1\leq i\leq n-1 , 1\leq j\leq m$
$$\lambda(u_i^ju_{i+1} ^ j)=v+e-m(n-1)-m(i-1)-j+1$$
$$\lambda(v_i^jv_{i+1} ^ j)=v+e-m(i-1)-j+1$$

\noindent
 Let $C^{(i,j)}_{4}:1\leq i \leq n-1$ and $1\leq j \leq m$, be the subcycle of
 graph $G$ with
 $V(C^{(i,j)}_{4})=\{u_i^j,u_{i+1}^j,v_i^j,v_{i+1}^j\}$ and $E(C^{(i,j)}_{4})=\{u_i^ju_{i+1}^j, v_i^jv_{i+1} ^ j,u_i^jv_i^j, u_{i+1} ^ jv_{i+1} ^ j\}.$
 After combining all three steps, we can verify that $\sum\ C^{(i,j)}_{4}=\lambda (u_i^j)+\lambda (u_{i+1}^j)+ \lambda
(v_i^j)+\lambda (v_{i+1}^j)+\lambda(u_i^jv_{i} ^
j)+\lambda(u_i^ju_{i+1} ^ j)+\lambda(v_i^jv_{i+1} ^
j)+\lambda(u_{i+1}^jv_{i+1} ^ j)=m(17n-2)+4.$  Hence $mL_n$ is $C_4$-supermagic.
 \qed\vspace{.5cm}

 \noindent Let $H\cong TL_n$ be a {\em triangular ladder} graph with
 $V(H)=\{u_i, v_i:\hphantom{a}1\leq i \leq
n\},$ $E(H)=\{u_iv_i:1\leq i \leq
n\}\cup\{u_iu_{i+1},v_iv_{i+1},u_{i+1}v_i:1\leq i \leq n-1\}$.

 \begin{theorem}
The graph $G\cong mTL_n$ admits
$C_3$-supermagic labelings for $m\geq 2$ and $n\geq 3.$
\end{theorem}
\noindent{\bf Proof.}\\
Let $v=|V(G)|$ and $e=|E(G)|$. Then $v=2mn$, $e=m(4n-3)$. We denote
the vertex and edge sets of
$G$ as follows:\\

 \noindent $V(G)=\{u^{j}_{i}: 1\leq i\leq n , 1\leq j\leq m \}\cup \{v^{j}_{i}: 1\leq i\leq n ,
1\leq j\leq m \},$\\
  $E(G)=\{u^{j}_{i}v^{j}_{i}: 1\leq i\leq n , 1\leq j\leq m \} \cup
\{u^{j}_{i}u^{j}_{i+1}:1\leq i\leq n-1 , 1\leq j\leq m\}\cup
\{v^{j}_{i}v^{j}_{i+1}:1\leq i\leq n-1 , 1\leq j\leq m\}\cup
\{u^{j}_{i+1}v^{j}_{i}: 1\leq i\leq n-1 , 1\leq j\leq m \}.$\\

\noindent Define a total labeling $\lambda:V(G)\cup E(G)\rightarrow
\{1,2,...,6mn-3m\}$ as follows:\\

\noindent For $1\leq i\leq n , 1\leq j\leq m.$
$$\lambda (u_i^j)=j+2m(i-1)$$
$$\lambda (v_i^j)=j+m(2i-1)$$
$$\lambda(u_i^jv_{i} ^ j)=m(4n+1-2i)-j+1$$

\noindent For $1\leq i\leq n-1 , 1\leq j\leq m$

$$\lambda(u_i^ju_{i+1} ^ j)=m(6n-3)-2m(i-1)-j+1$$
$$\lambda(v_i^jv_{i+1} ^ j)=2m(3n-i-1)-j+1$$
$$\lambda(u_{i+1}^jv_{i} ^ j)=2m(2n-i)-j+1$$

\noindent It is easy to verify that for every subcycle $C^{(i,j)}_{3}:
1\leq i \leq 2n-2$  and $1\leq j \leq m$, $\lambda
(u_i^j)+\lambda (u_{i+1}^j)+ \lambda (v_i^j)+\lambda(u_i^jv_{i} ^
j)+\lambda(u_i^ju_{i+1} ^ j)+\lambda(u_{i+1}^jv_{i} ^ j)=14mn-3m+3.$
 Hence $mTL_n$ is $C_3$-supermagic.\qed\vspace{.5cm}

\noindent Next we consider $C_3$-supermagic labeling of disjoint union of isomorphic copies of odd wheel graphs.
 \noindent
 Let $W_n =K_1+C_n $ be a wheel graph with
 $V(W_n)=\{c\} \cup \{ v_i:\hphantom{a}1\leq i \leq
n\},$ and
 $E(W_n)=\{cv_i:1\leq i \leq n\}\cup\{v_i v_{i+1},1\leq i \leq
n\}.$
\\

\begin{theorem}
The graph
$G\cong mW_n$ is $C_3$-supermagic, for any odd $n\geq 3$ and positive integer $m\geq 2.$
\end{theorem}
\noindent{\bf Proof.}\\
Let $v=|V(G)|$ and $e=|E(G)|$. Then $v=m(n+1)$, $e=2mn$. We denote
the vertex and edge sets of $G$ as follows:\\

\noindent $V(G)=\{c_{j};1\leq j \leq m \}\cup \{v^{j}_{i}: 1\leq i\leq n ,
1\leq j\leq m \},$ and\\
\noindent $E(G)=\{v^{j}_{i} v^{j}_{i+1}: 1\leq i\leq n-1 , 1\leq j\leq m \} \cup
\{c_j v^{j}_{i}:1\leq i\leq n , 1\leq j\leq m \}  .$\\

\noindent
 Define a total labeling $\lambda:V(G)\cup E(G)\rightarrow
\{1,2,...,m(3n+1)\}$ as follows:\\

\noindent For $1\leq i\leq n , 1\leq j\leq m$
$$\lambda (c_j)=j$$
$$\lambda (c_j v_i^j)=m(2n-i+1)-j+1$$
$$\lambda (c_j v_n)=m(2n+1)-j+1$$
$$\lambda (v_nv_1)=m(2n+3)-j+1$$
$$\lambda(v_i^j v_{i+1} ^ j)=m(2n+i+3)-j+1$$
$$
\lambda (v_i^j)= \left\{ \begin{array}{l@{\quad}l} m(\frac{i+1}{2}) +j,
&\hphantom{i} \mbox{}\hphantom{i} $if$ \hphantom{i} i -odd\\\\
\mbox{ } m(\frac{n+i+1}{2})+j, & \mbox{}
$if$ \hphantom{i} i-even\\
\end{array} \right.$$

\noindent It is easy to verify that for every subcycle $C^{(i,j)}_{3}: 1\leq i
\leq n,\,1\leq j\leq m$ of $mW_n$, the weight of every subcycle $C^{(i,j)}_{3}$ of $G$ is $\frac{m}{2}(13n+11)+3$.
Hence $mW_n$ is $C_{3}$-supermagic .\qed\vspace{.5cm}

\begin{theorem}
 For even $n\geq 4,$ and
$m\geq 2,$ the graph $G\cong mW_n$ is $C_3$-supermagic.
\end{theorem}
\noindent{\bf Proof.}\\
Let $v=|V(G)|$ and $e=|E(G)|$. Then $v=m(n+1)$, $e=2mn$. We denote
the vertex and edge sets of $G$ as follows:\\

\noindent $V(G)=\{v_0^j:\,1\leq j \leq m \}\cup \{v^{j}_{i}: \,1\leq
i\leq n ,
1\leq j\leq m \},$ and\\
\noindent $E(G)=\{v^{j}_{i} v^{j}_{i+1}:\, 1\leq i\leq n-1 , 1\leq
j\leq m \} \cup \{v_0^j v^{j}_{i}:\,1\leq i\leq n , 1\leq j\leq m
\}\cup\{v_n^j v^{j}_{1}:\, 1\leq j\leq m \}  .$\\

\noindent
 Define a total labeling $\lambda:V(G)\cup E(G)\rightarrow
\{1,2,...,m(3n+1)\}$ as follows:\\
\noindent Throughout the following labeling we will consider $1\leq
j\leq m$,

\noindent{\bf Case 1: $\frac{n}{2}$-even}
$$\lambda (v_0^j)=m(\frac{n}{4}+\lceil \frac{n-1}{2} \rceil+1)+j-m$$
$$
\lambda (v_i^j)= \left\{ \begin{array}{l@{\quad}l} m\lfloor
\frac{i}{2}\rfloor+j, &\hphantom{i}
\mbox{}\hphantom{i} $if$ \hphantom{i} i-odd, 1\leq i\leq n\\\\
\frac{m}{2}(i+n)+j-m, &\hphantom{i}
\mbox{}\hphantom{i} $if$ \hphantom{i} i-even, 1\leq i\leq \frac{n}{2}\\\\
\frac{m}{2}(i+n)+j, &\hphantom{i}
\mbox{}\hphantom{i} $if$ \hphantom{i} i-even, \frac{n}{2}+1 \leq i \leq n\\\\
\end{array} \right.$$
$$
\lambda (v_0^j v_i^j)= \left\{ \begin{array}{l@{\quad}l}
m(2n+2-i)-(j-1), &\hphantom{i}
\mbox{}\hphantom{i} $if$ \hphantom{i} 1\leq i\leq \frac{n}{2}\\\\
m(2n+1-i)-(j-1), &\hphantom{i}
\mbox{}\hphantom{i} $if$ \hphantom{i} \frac{n}{2}+1 \leq i \leq n-1\\\\
m(\frac{3n}{2}+1)-(j-1), &\hphantom{i}
\mbox{}\hphantom{i} $if$ \hphantom{i} i=n\\\\
\end{array} \right.$$
$$
\lambda (v_i^j v_{i+1}^j)= \left\{ \begin{array}{l@{\quad}l}
m(2n+2+i)-(j-1), &\hphantom{i}
\mbox{}\hphantom{i} $if$ \hphantom{i} 1\leq i\leq \frac{n}{2}-1\\\\
m(2n+3+i)-(j-1), &\hphantom{i}
\mbox{}\hphantom{i} $if$ \hphantom{i} \frac{n}{2} \leq i \leq n-2\\\\
m(\frac{5n}{2}+2)-(j-1), &\hphantom{i}
\mbox{}\hphantom{i} $if$ \hphantom{i} i=n-1\\\\
\end{array} \right.$$
$$\lambda (v_n^j v_1^j)=2m(n+1)-(j-1)$$

\noindent It is easy to verify that weight of every subcycle
$C^{(i,j)}_{3}: 1\leq i \leq n,\,1\leq j\leq m$ of $mW_n$ is $m
\lceil \frac{n-1}{2}\rceil +\frac{m}{4}(27n+16)+3$. Hence $mW_n$ is
$C_{3}$-supermagic .\\

\noindent{\bf Case 2: $\frac{n}{2}$-odd}
$$\lambda (v_0^j)=\frac{3m(n+2)}{4}+j-m$$
$$
\lambda (v_i^j)= \left\{ \begin{array}{l@{\quad}l} m\lfloor
\frac{i}{2}\rfloor +j, &\hphantom{i}
\mbox{}\hphantom{i} $if$ \hphantom{i} i-odd, 1\leq i\leq n\\\\
m[\frac{1}{2}(i+n)+1]+j-m, &\hphantom{i}
\mbox{}\hphantom{i} $if$ \hphantom{i} i-even, 1\leq i\leq \frac{n}{2}\\\\
m[\frac{1}{2}(i+n)+2]+j-m, &\hphantom{i}
\mbox{}\hphantom{i} $if$ \hphantom{i} i-even, \frac{n}{2}+1 \leq i \leq n-1\\\\
m(\frac{i}{2})+j, &\hphantom{i}
\mbox{}\hphantom{i} $if$ \hphantom{i} i=n\\\\
\end{array} \right.$$
$$
\lambda (v_0^j v_i^j)= \left\{ \begin{array}{l@{\quad}l}
2m(n+1)-(j-1), &\hphantom{i}
\mbox{}\hphantom{i} $if$ \hphantom{i} i=1\\\\
m(2n+2-i)-(j-1), &\hphantom{i} \mbox{}\hphantom{i} $if$ \hphantom{i}
2\leq i\leq \frac{n}{2}\\\\
m(2n-i)-(j-1), &\hphantom{i}
\mbox{}\hphantom{i} $if$ \hphantom{i} i-even, \frac{n}{2}+1 \leq i \leq n-1\\\\
m(2n+2-i)-(j-1), &\hphantom{i}
\mbox{}\hphantom{i} $if$ \hphantom{i} i-odd, \frac{n}{2}+1 \leq i \leq n\\\\
m(\frac{3n}{2}+1)-(j-1), &\hphantom{i}
\mbox{}\hphantom{i} $if$ \hphantom{i} i=n\\\\
\end{array} \right.$$

$$
\lambda (v_i^j v_{i+1}^j)= \left\{ \begin{array}{l@{\quad}l}
m(2n+1)-(j-1), &\hphantom{i}
\mbox{}\hphantom{i} $if$ \hphantom{i} i=1\\\\
m(2n+1+i)-(j-1), &\hphantom{i}
\mbox{}\hphantom{i} $if$ \hphantom{i} 2 \leq i \leq \frac{n}{2}-1\\\\
m(2n+2+i)-(j-1), &\hphantom{i}
\mbox{}\hphantom{i} $if$ \hphantom{i} \frac{n}{2} \leq i \leq n-1\\\\
\end{array} \right.$$
$$\lambda (v_n^j v_1^j)=m(\frac{5n}{2}+1)-(j-1)$$

\noindent The weight of every subcycle $C^{(i,j)}_{3}: 1\leq i \leq
n,\,1\leq j\leq m$ of $mW_n$, under the labeling $\lambda$  is
$\frac{m}{4}(29n+18)+3$. Hence $mW_n$ is
$C_{3}$-supermagic.\qed\vspace{.5cm}

\noindent The book graph $B_n \cong K_{1,n}\times K_2$ has the
vertex set $V(G)=\{u_1, u_2, v_i, w_i \}$ and the edge set
$E(G)=\{u_1, u_2, u_1 w_i, u_2 v_i, v_i w_i \}$ for $1\leq i\leq
n$.\\\\

\noindent In the following theorem, we consider $C_4$-supermagic
labeling for disjoint union of isomorphic copies of book
graphs.\\\\

 \begin{theorem}
For $m\geq2$, $n\geq 2$ the graph  $G\cong mB_n$ admits
$C_4$-supermagic labelings.
\end{theorem}
\noindent{\bf Proof.}\\
Let $G$ be a graph and let $v=|V(G)|$ and $e=|E(G)|$. Then
$v=2m(n+1)$, $e=m(3n+1)$. We denote the vertex and edge sets of $G$
as follows:\\

\noindent $V(G)=\{u^{j}_{1}, u^{j}_{2}: 1\leq j\leq m\} \cup \{
v^{j}_{i}, w^{j}_{i} : 1\leq i\leq n , 1\leq j\leq m \}$\\
$E(G)=\{u^{j}_{1} u^{j}_{2}: 1\leq j\leq m \} \cup
\{u^{j}_{1}w^{j}_{i}, u^{j}_{2}v^{j}_{i}:1\leq i\leq n , 1\leq j\leq
m \}\cup \{ v^{j}_{i}w^{j}_{i}: 1\leq i\leq n , 1\leq j\leq m \}.$

\noindent Throughout the following labeling, we will consider $1\leq
i\leq n , 1\leq j\leq m.$
$$\lambda (u^{j}_{i})= m(i-1)+j,  i=1,2$$
$$\lambda (v^{j}_{i})= m(i+1)+j$$
$$\lambda (w^{j}_{i})= m(2n+2-i)+j.$$
\noindent Further, we will split the labelings into two cases.\\
\noindent { \bf Case 1: n is even.}
$$
\lambda (u^{j}_{1} u^{j}_{2})= m(\frac{5n}{2}+3)-j+1
$$
$$
\lambda (u^{j}_{2} v^{j}_{i})= \left\{ \begin{array}{l@{\quad}l}
m(2n+2+i)-j+1, &\hphantom{i}
\mbox{}\hphantom{i} $if$ \hphantom{i}  1\leq i\leq \frac{n}{2}\\\\
m(2n+3+i)-j+1, &\hphantom{i}
\mbox{}\hphantom{i} $if$ \hphantom{i}  \frac{n}{2}+1 \leq i\leq n\\\\
\end{array} \right.$$
$$
\lambda (u^{j}_{1} w^{j}_{i})= \left\{ \begin{array}{l@{\quad}l}
m(5n+5-2i)-j+1, &\hphantom{i}
\mbox{}\hphantom{i} $if$ \hphantom{i}  1\leq i\leq \frac{n}{2}\\\\
m(6n+4-2i)-j+1, &\hphantom{i}
\mbox{}\hphantom{i} $if$ \hphantom{i}  \frac{n}{2}+1 \leq i\leq n\\\\
\end{array} \right.$$
$$
\lambda (v^{j}_{i} w^{j}_{i})= \left\{ \begin{array}{l@{\quad}l}
(\frac{m}{2})(7n+6+2i)-j+1, &\hphantom{i}
\mbox{}\hphantom{i} $if$ \hphantom{i}  1\leq i\leq \frac{n}{2}\\\\
(\frac{m}{2})(5n+6+2i)-j+1, &\hphantom{i}
\mbox{}\hphantom{i} $if$ \hphantom{i}  \frac{n}{2}+1 \leq i\leq n\\\\
\end{array} \right.$$

\noindent It is easy to verify that for every subcycle $C^{(i)}_{4}:
1\leq i \leq n$  and $1\leq j \leq m$, the weight is $\lambda
(u_1^j)+\lambda (u_2^j)+ \lambda (v_i^j)+\lambda
(w_i^j)+\lambda(u_1^ju_2^ j)+\lambda(u_2^jv_i^ j)+\lambda(u_1^jw_i ^
j)+\lambda(v_i^jw_{i} ^ j)=15mn+17m+4$. Hence $mB_n$ for $n$-even is
$C_4$-supermagic.\\

 \noindent {\bf Case 2: n is odd.}
$$
\lambda (u^{j}_{1} u^{j}_{2})= m(2n+3)-j+1
$$
$$
\lambda (u^{j}_{2} v^{j}_{i})= m(2n+3+i)-j+1
$$
$$
\lambda (u^{j}_{1} w^{j}_{i})= \left\{ \begin{array}{l@{\quad}l}
m(5n+5-2i)-j+1, &\hphantom{i}
\mbox{}\hphantom{i} $if$ \hphantom{i}  1\leq i\leq \lceil\frac{n}{2}\rceil\\\\
m(6n+5-2i)-j+1, &\hphantom{i}
\mbox{}\hphantom{i} $if$ \hphantom{i}  \lceil\frac{n}{2}\rceil+1 \leq i\leq n\\\\
\end{array} \right.$$
$$
\lambda (v^{j}_{i} w^{j}_{i})= \left\{ \begin{array}{l@{\quad}l}
(\frac{m}{2})(7n+5+2i)-j+1, &\hphantom{i}
\mbox{}\hphantom{i} $if$ \hphantom{i}  1\leq i\leq \lceil\frac{n}{2}\rceil\\\\
(\frac{m}{2})(5n+5+2i)-j+1, &\hphantom{i}
\mbox{}\hphantom{i} $if$ \hphantom{i}  \lceil\frac{n}{2}\rceil+1 \leq i\leq n\\\\
\end{array} \right.$$

\noindent It is easy to verify that for every subcycle $C^{(i)}_{4}:
1\leq i \leq n$, $1\leq j \leq m$, the weigh is $\lambda
(u_1^j)+\lambda (u_2^j)+ \lambda (v_i^j)+\lambda
(w_i^j)+\lambda(u_1^ju_2^ j)+\lambda(u_2^jv_i^ j)+\lambda(u_1^jw_i ^
j)+\lambda(v_i^jw_{i} ^ j)=\frac{m}{2}(29n+35)+4$ . Hence $mB_n$ for
$n$-odd is $C_4$-supermagic. \qed\vspace{.5cm}

\noindent The prism is a graph $G\cong C_m\times P_n$ with
vertex-set $V(G) =\{v_{i,j}:1\leq i \leq m,1\leq j \leq n \} $ and
edge-set $E(G)=\{v_{i,j}v_{i,j+1}: 1\leq i\leq m, 1\leq j\leq n-1
\}\cup
\{v_{i,j}v_{i+1,j},1\leq i\leq m , 1\leq j\leq n\}$.\\
\noindent A generalized antiprism $A^{n}_{m}$ is a graph obtained by
completing the prism graph $C_m\times P_n$ by adding the edges
$v_{i,j+1}v_{i+1,j}$ for $1\leq i\leq m$ and $1\leq j\leq n-1.$

\begin{theorem}
For $l\geq 2$ and $m, n\geq 3 $, the generalized antiprism $
A^{n}_{m}$ is $C_3$-supermagic with magic constant $lm(9n-4)+3$.
\end{theorem}

\noindent{\bf Proof.}\\

\noindent Let $v=|V(A^{n}_{m})|$ and $e=|E(A^{n}_{m})|$. Then
$v=lmn$, $e=lm(3n-2)$. We denote the vertex and edge-sets of
$G$ as follows:\\

 \noindent $V(A^{n}_{m})=\{v^{k}_{i,j}: 1\leq i\leq m , 1\leq j\leq n, 1\leq k\leq l \},$\\
  $E(A^{n}_{m})=\{v^{k}_{i,j}v^{k}_{i,j+1}: 1\leq i\leq m , 1\leq j\leq n-1, 1\leq k\leq l \} \cup
\{v^{k}_{i,j+1}v^{k}_{i+1,j},1\leq i\leq m ,1\leq j\leq n-1, 1\leq
k\leq l\}\cup \{v^{k}_{i,j}v^{k}_{i+1,j},1\leq i\leq m , 1\leq j\leq
n, 1\leq k\leq l\}.$\\

\noindent Let $B^k_{i,j}$ be the $3$-cycle
$v_{i,j}^kv_{i+1,j}^kv_{i,j+1}^k$ and $C^{k}_{i,j}$ be the $3$-cycle
$v_{i,j+1}^kv_{i+1,j+1}^kv_{i+1,j}^k$.\\

\noindent Define a total labeling $\lambda:V(A^{n}_{m})\cup
E(A^{n}_{m})\rightarrow
\{1,2,...,4lmn-2lm\}$ as follows:\\
\noindent{\bf Case 1: $m$-odd}

$$
\lambda (v_{i,j}^k)= \left\{ \begin{array}{l@{\quad}l}
l[m(j-1)+\frac{m-i+2}{2}]+k-l, &\hphantom{i}
\mbox{}\hphantom{i} $if$ \hphantom{i} i-odd, j-odd\\\\
l[m(j-1)+\frac{2m-i+2}{2}]+k-l, &\hphantom{i}
\mbox{}\hphantom{i} $if$ \hphantom{i} i-even, j-odd\\\\
l[m(j-1)+\frac{i+1}{2}]-k+1, &\hphantom{i}
\mbox{}\hphantom{i} $if$ \hphantom{i} i-odd, j-even\\\\
l[m(j-1)+\frac{m+i+1}{2}]-k+1, &\hphantom{i}
\mbox{}\hphantom{i} $if$ \hphantom{i} i-even, j-even\\\\
\end{array} \right.$$
$$
\lambda (v_{i,j}^k v_{i+1,j}^k)= \left\{ \begin{array}{l@{\quad}l}
l[m(4n-j-2)+\frac{i+1}{2}]-k+1, &\hphantom{i}
\mbox{}\hphantom{i} $if$ \hphantom{i} i-odd, j-odd\\\\
l[m(4n-j-2)+\frac{m+i+1}{2}]-k+1, &\hphantom{i}
\mbox{}\hphantom{i} $if$ \hphantom{i} i-even, j-odd\\\\
l[m(4n-j-2)+\frac{m-i}{2}]+k-l, &\hphantom{i}
\mbox{}\hphantom{i} $if$ \hphantom{i} i-odd, j-even, i\neq m\\\\
lm(4n-j-1)+k-l, &\hphantom{i}
\mbox{}\hphantom{i} $if$ \hphantom{i} i=m, j-even\\\\
l[m(4n-j-2)+\frac{2m-i}{2}]+k-l, &\hphantom{i}
\mbox{}\hphantom{i} $if$ \hphantom{i} i-even, j-even\\\\
\end{array} \right.$$

\noindent For $n$-even
$$
\lambda (v_{i,j}^k v_{i,j+1}^k)= \left\{ \begin{array}{l@{\quad}l}
l[m(2n-j-1)+\frac{m+i}{2}]+k-l, &\hphantom{i}
\mbox{}\hphantom{i} $if$ \hphantom{i} i-odd, j-odd\\\\
l[m(2n-j-1)+\frac{i}{2}]+k-l, &\hphantom{i}
\mbox{}\hphantom{i} $if$ \hphantom{i} i-even, j-odd\\\\
l[m(2n-j-1)+\frac{2m-i+1}{2}]+k-l, &\hphantom{i}
\mbox{}\hphantom{i} $if$ \hphantom{i} i-odd, j-even\\\\
l[m(2n-j-1)+\frac{m-i+1}{2}]+k-l, &\hphantom{i}
\mbox{}\hphantom{i} $if$ \hphantom{i} i-even, j-even\\\\
\end{array} \right.$$

\noindent For $n$-odd
$$
\lambda (v_{i,j}^k v_{i,j+1}^k)= \left\{ \begin{array}{l@{\quad}l}
l[m(2n-j-1)+\frac{m+i}{2}]-k+1, &\hphantom{i}
\mbox{}\hphantom{i} $if$ \hphantom{i} i-odd, j-odd\\\\
l[m(2n-j-1)+\frac{i}{2}]-k+1, &\hphantom{i}
\mbox{}\hphantom{i} $if$ \hphantom{i} i-even, j-odd\\\\
l[m(2n-j-1)+\frac{2m-i+1}{2}]-k+1, &\hphantom{i}
\mbox{}\hphantom{i} $if$ \hphantom{i} i-odd, j-even\\\\
l[m(2n-j-1)+\frac{m-i+1}{2}]-k+1, &\hphantom{i}
\mbox{}\hphantom{i} $if$ \hphantom{i} i-even, j-even\\\\
\end{array} \right.$$

\noindent For $n$-even
$$
\lambda (v_{i,j+1}^k v_{i+1,j}^k)= \left\{
\begin{array}{l@{\quad}l} l[m(3n-j-2)+\frac{2m-i+1}{2}]-k+1,
&\hphantom{i}
\mbox{}\hphantom{i} $if$ \hphantom{i} i-odd, j-odd\\\\
l[m(3n-j-2)+\frac{m-i+1}{2}]-k+1, &\hphantom{i}
\mbox{}\hphantom{i} $if$ \hphantom{i} i-even, j-odd\\\\
l[m(3n-j-2)+\frac{m+i}{2}]-k+1, &\hphantom{i}
\mbox{}\hphantom{i} $if$ \hphantom{i} i-odd, j-even\\\\
l[m(3n-j-2)+\frac{i}{2}]-k+1, &\hphantom{i}
\mbox{}\hphantom{i} $if$ \hphantom{i} i-even, j-even\\\\
\end{array} \right.$$

\noindent For $n$-odd
$$
\lambda (v_{i,j+1}^k v_{i+1,j}^k)= \left\{
\begin{array}{l@{\quad}l} l[m(3n-j-2)+\frac{2m-i+1}{2}]+k-l,
&\hphantom{i}
\mbox{}\hphantom{i} $if$ \hphantom{i} i-odd, j-odd\\\\
l[m(3n-j-2)+\frac{m-i+1}{2}]+k-l, &\hphantom{i}
\mbox{}\hphantom{i} $if$ \hphantom{i} i-even, j-odd\\\\
l[m(3n-j-2)+\frac{m+i}{2}]+k-l, &\hphantom{i}
\mbox{}\hphantom{i} $if$ \hphantom{i} i-odd, j-even\\\\
l[m(3n-j-2)+\frac{i}{2}]+k-l, &\hphantom{i}
\mbox{}\hphantom{i} $if$ \hphantom{i} i-even, j-even\\\\
\end{array} \right.$$

\noindent For $1\leq i\leq m$, $1\leq j\leq n-1$ and $1\leq k\leq
l$,
$$ \sum (B_{i,j}^k)= \lambda (v_{i,j}^k)+ \lambda (v_{i+1,j}^k)+\lambda
(v_{i,j+1}^k)+$$
$$+\lambda (v_{i,j}^k v_{i+1,j}^k)+\lambda (v_{i,j+1}^k
v_{i+1,j}^k)+\lambda (v_{i,j}^k v_{i,j+1}^k)
$$
where $i$ is odd and $j$ is odd,
$$ \sum(B_{i,j}^k)= l[m(j-1)+\frac{m-i+2}{2}]+k-l+l[m(j-1)+$$
$$+\frac{2m-i+2}{2}]+k-l+ l[mj+\frac{i+1}{2}]-k+1+$$
$$+l[m(4n-j-2)+\frac{i+1}{2}]-k+1+ l[m(2n-j-1)+$$
$$+\frac{m+i}{2}]-k+1+l[m(3n-j-2)+\frac{2m-i+1}{2}]+k-l$$
$$\sum (B_{i,j}^k)=lm[9n-4]+3$$
where $i$ is even and $j$ is even,
$$ \sum (B_{i,j}^k)= l[m(j-1)+\frac{m+i+1}{2}]-k+1+l[m(j-1)+$$
$$+\frac{i+2}{2}]-k+1+ l[m(j-1)+\frac{2m-i+2}{2}]+k-l+$$
$$+l[m(4n-j-2)+\frac{2m-i}{2}]+k-l+ l[m(2n-j-1)+$$
$$+\frac{m-i+1}{2}]+k-l+l[m(3n-j-2)+\frac{i}{2}]-k+1$$
$$\sum (B_{i,j}^k)=lm[9n-4]+3$$

\noindent Thus, it can be easily verified that $ \sum
(B_{i,j}^k)=lm(9n-4)+3 $ for the remaining cases of $i$ and $j$,
where $1\leq i\leq m, 1\leq j\leq n-1$ and $1\leq k\leq l.$\\

\noindent Next we prove $ \sum (C_{i,j}^k)=lm(9n-4)+3 $, where
$1\leq i\leq m , 1\leq j\leq n-1$ and $1\leq k\leq l.$

$$ \sum (C_{i,j}^k)= \lambda (v_{i,j+1}^k)+ \lambda (v_{i+1,j}^k)+\lambda
(v_{i+1,j+1}^k)+$$
$$+\lambda (v_{i,j+1}^k v_{i+1,j}^k)+\lambda (v_{i+1,j+1}^k
v_{i+1,j+1}^k)+\lambda (v_{i,j+1}^k v_{i+1,j+1}^k)
$$
where $i$ is odd and $j$ is odd and $i\neq m,$
$$ \sum (C_{i,j}^k)= l[mj+\frac{m+i+1}{2}]-k+1+l[m(j-1)+$$
$$+\frac{2m-i+1}{2}]+k-l+ l[mj+\frac{m+i+2}{2}]-k+1+$$
$$+l[m(4n-j-3)+\frac{m-i}{2}]+k-l+ l[m(2n-j-1)+$$
$$+\frac{i+1}{2}]-k+1+l[m(3n-j-2)+\frac{2m-i+1}{2}]+k-l$$
$$\sum (C_{i,j}^k)=lm[9n-4]+3$$
where $i$ is even and $j$ is even,

$$ \sum (C_{i,j}^k)= l[mj+\frac{2m-i+2}{2}]+k-l+l[m(j-1)+$$
$$+\frac{i+2}{2}]-k+1+ l[mj+\frac{m-i+1}{2}]+k-l+$$
$$+l[m(4n-j-3)+\frac{m+i+1}{2}]-k+1+ l[m(2n-j-1)+$$
$$+\frac{2m-i}{2}]+k-l+l[m(3n-j-2)+\frac{i}{2}]-k+1$$
$$\sum (C_{i,j}^k)=lm[9n-4]+3$$

\noindent Similarly for the other cases of $i$ and $j$, $\sum
(B_{i,j}^k)=\sum (C_{i,j}^k)=lm[9n-4]+3$, where $1\leq i\leq m,
1\leq j\leq n-1$ and $1\leq k\leq l.$ Thus we get a $C_3$-supermagic
labeling of $A_m^n$ with supermagic constant $c=lm[9n-4]+3$.
Hence $A^{n}_{m}$ is $C_3$-supermagic.\\

\noindent{\bf Case 2: $m$-even}

$$
\lambda (v_{i,j}^k)= \left\{ \begin{array}{l@{\quad}l}
l[m(j-1)+i]+k-l, &\hphantom{i}
\mbox{}\hphantom{i} $if$ \hphantom{i}  j-odd\\\\
l[mj-i+1)]-k+1, &\hphantom{i}
\mbox{}\hphantom{i} $if$ \hphantom{i}  j-even\\\\
\end{array} \right.$$
$$
\lambda (v_{i,j}^k v_{i+1,j}^k)= \left\{ \begin{array}{l@{\quad}l}
l[m(4n-j-1)-i]-k+1, &\hphantom{i}
\mbox{}\hphantom{i} $if$ \hphantom{i}   j-odd, i\neq m\\\\
lm(4n-j-1)-k+1, &\hphantom{i}
\mbox{}\hphantom{i} $if$ \hphantom{i}   j-odd, i= m\\\\
l[m(4n-j-2)+i+1]+k-l, &\hphantom{i}
\mbox{}\hphantom{i} $if$ \hphantom{i}  j-even, i\neq m\\\\
l[m(4n-j-2)+1]+k-l, &\hphantom{i}
\mbox{}\hphantom{i} $if$ \hphantom{i}  j-even, i=m\\\\
\end{array} \right.$$

$$
\lambda (v_{i,j}^k v_{i,j+1}^k)= \left\{ \begin{array}{l@{\quad}l}
l[m(2n-j)-i+1]+k-l, &\hphantom{i}
\mbox{}\hphantom{i} $if$ \hphantom{i}  j-odd, n-even\\\\
l[m(2n-j)-i+1]-k+1, &\hphantom{i}
\mbox{}\hphantom{i} $if$ \hphantom{i}  j-odd, n-odd\\\\
l[m(2n-j-1)+i]-k+1, &\hphantom{i}
\mbox{}\hphantom{i} $if$ \hphantom{i}  j-even, n-odd\\\\
l[m(2n-j-1)+i]+k-l, &\hphantom{i}
\mbox{}\hphantom{i} $if$ \hphantom{i}  j-even,  n-even\\\\
\end{array} \right.$$
$$
\lambda (v_{i,j+1}^k v_{i+1,j}^k)= \left\{
\begin{array}{l@{\quad}l} l[m(3n-j-2)+i]-k+1, &\hphantom{i}
\mbox{}\hphantom{i} $if$ \hphantom{i}  j-odd, n-even\\\\
l[m(3n-j-2)+i]+k-l, &\hphantom{i}
\mbox{}\hphantom{i} $if$ \hphantom{i}  j-odd, n-odd\\\\
l[m(3n-j-1)-i+1]-k+1, &\hphantom{i}
\mbox{}\hphantom{i} $if$ \hphantom{i}  j-even, n-even\\\\
l[m(3n-j-1)-i+1]+k-l, &\hphantom{i}
\mbox{}\hphantom{i} $if$ \hphantom{i}  j-even, n-odd\\\\
\end{array} \right.$$\\
\noindent For $1\leq i\leq m$, $1\leq j\leq n-1$ and $1\leq k\leq
l,$

$$ \sum (B_{i,j}^k)= \lambda (v_{i,j}^k)+ \lambda (v_{i+1,j}^k)+\lambda
(v_{i,j+1}^k)+$$
$$+\lambda (v_{i,j}^k v_{i+1,j}^k)+\lambda (v_{i,j+1}^k
v_{i+1,j}^k)+\lambda (v_{i,j}^k v_{i,j+1}^k)
$$
where $j$ is odd and $i \neq m$,
$$ \sum (B_{i,j}^k)= l[m(j-1)+i]+k-l+l[m(j-1)+i+1]+k-l$$
$$+l[m(j+1)-i+1]-k+1+l[m(4n-j-1)-i]-k+1+$$
$$+l[m(2n-j)-i+1]+k-l+l[m(3n-j-2)+i]-k+1$$
$$\sum (B_{i,j}^k)=lm[9n-4]+3$$

where $j$ is even and $i=m$,
$$ \sum (B_{m,j}^k)= \lambda (v_{m,j}^k)+ \lambda (v_{1,j}^k)+\lambda
(v_{m,j+1}^k)+$$
$$+\lambda (v_{m,j}^k v_{m+1,j}^k)+\lambda (v_{m,j+1}^k
v_{m+1,j}^k)+\lambda (v_{m,j}^k v_{m,j+1}^k)
$$
$$ \sum (B_{i,j}^k)= l[mj-m+1]-k+1+lmj-k+1$$
$$+l[mj+m]+k-l+l[m(4n-j-2)+1]+k-l+$$
$$+l[m(2n-j-1)+m]+k-l+l[m(3n-j-1)-m+1]-k+1$$
$$\sum (B_{i,j}^k)=lm[9n-4]+3$$

\noindent Similarly for the remaining cases of $i$ and $j$, it can
be easily verified that $ \sum (B_{i,j}^k)=lm(9n-4)+3 $ for $1\leq
i\leq m, 1\leq j\leq n-1$ and $1\leq k\leq l.$

\noindent Next we prove $ \sum (C_{i,j}^k)=lm(9n-4)+3 $, where
$1\leq i\leq m , 1\leq j\leq n-1$ and $1\leq k\leq l.$
$$ \sum (C_{i,j}^k)= \lambda (v_{i,j+1}^k)+ \lambda (v_{i+1,j}^k)+\lambda
(v_{i+1,j+1}^k)+$$
$$+\lambda (v_{i,j+1}^k v_{i+1,j}^k)+\lambda (v_{i+1,j}^k
v_{i+1,j+1}^k)+\lambda (v_{i,j+1}^k v_{i+1,j+1}^k)
$$
where $j$ is odd and $i \neq m$,
$$ \sum (C_{i,j}^k)= l[m(j+1)-i+1]-k+1+l[m(j-1)+i+1]+k-l$$
$$+l[m(j+1)-i]-k+1+l[m(3n-j-2)+i]-k+1+$$
$$+l[m(2n-j)-i]+k-l+l[m(4n-j-3)+i+1]+k-l$$
$$\lambda (C_{i,j}^k)=lm[9n-4]+3$$

where $j$ is even and $i=m$,
$$ \sum (C_{m,j}^k)= \lambda (v_{m,j+1}^k)+ \lambda (v_{1,j}^k)+\lambda
(v_{1,j+1}^k)+$$
$$+\lambda (v_{m,j+1}^k v_{m+1,j}^k)+\lambda (v_{1,j}^k
v_{1,j+1}^k)+\lambda (v_{m,j+1}^k v_{m+1,j+1}^k)
$$
$$ \sum (C_{m,j}^k)= l[mj+m]+k-l+lmj-k+1$$
$$+l[mj+1]+k-l+l[m(3n-j-1)+1-m]+k-l+$$
$$+l[m(2n-j-1)+1]-k+1+l[m(4n-j-2)-m+1]-k+1$$
$$\sum (C_{m,j}^k)=lm[9n-4]+3$$

\noindent Similarly for the other cases of $i$ and $j$, $\sum
(B_{i,j}^k)=\sum (C_{i,j}^k)=lm[9n-4]+3$, where $1\leq i\leq m ,
1\leq j\leq n-1$ and $1\leq k\leq l.$ Thus we get a $C_3$-supermagic
labeling of $A_m^n$ with supermagic constant $c=lm[9n-4]+3$. Hence
$A^{n}_{m}$ is $C_3$-supermagic. \qed\vspace{.5cm}

 \noindent In the
following section,  we study the cycle-supermagic labeling
 for the disjoint union of non isomorphic copies of   fans and ladders.
\section{Cycle-supermagic labelings of the disjoint union of non isomorphic graphs}
Let $G\cong sF_{n+1} \cup kF_n,$ $s,k\geq1,$ we denote the vertex and edge sets of $G$ as
follows:\\

 \noindent $V(G)=\{c_j: 1\leq j\leq s+k \}\cup \{v^{j}_{i}: 1\leq i\leq b , 1\leq j\leq s+k \},$\\
\noindent$E(G)=\{v^{j}_{i}v^{j}_{i+1}: 1\leq i\leq b , 1\leq j\leq s+k \}
\cup \{c_jv^{j}_{i}:1\leq i\leq b , 1\leq j\leq s+k \}, $
where $$
b= \left\{ \begin{array}{l@{\quad}l}
n &\hphantom{i} \mbox{}\hphantom{i} $if$
\hphantom{i} 1\leq j\leq s,\\\\ \mbox{ } n-1 &
\hphantom{i}\hphantom{i}\mbox{}
$if$ \hphantom{i} s+1\leq j\leq s+k.\\
\end{array} \right.$$

\begin{theorem}
For any positive integer $s,$ $k$ and $n\geq 3$, the graph
$G\cong sF_{n+1} \cup kF_n$ is $C_3$-supermagic.
\end{theorem}
\noindent{\bf Proof.}\\
Let $v=|V(G)|$ and $e=|E(G)|$. Then $v=s(n+1)+nk$,
$e=s(2n-1)+k(2n-3)$.

\noindent
 Define a total labeling $\lambda:V(G)\cup E(G)\rightarrow
\{1,2,....,n(3s+2k)-3k \}$ as follows:\\

\noindent Throughout the following labeling, we will consider $1\leq i\leq b, 1\leq
j\leq s+k.$
$$\lambda(c_j)=s(n+1)+nk-j+1$$
\noindent For $i$-odd

$$\lambda (v_i^j)=j+(\frac{s+k}{2})(i-1)$$

\noindent For $i$-even
$$\lambda (v_i^j)=j+(\frac{s+k}{2})(i-2)+\lfloor(\frac{s+k}{2})n\rfloor,$$

$$\lambda(c_jv_{i}^j)=s(3n-i+1)+k(3n-i-2)-j+1,$$
$$\lambda(v_i^jv_{i+1} ^ j)=s(n+i)+k(n+i-1)+j,\,\,\, 1\leq i\leq b-1 .$$

\noindent It is easy to verify that for every subcycle $C^{(i)}_{3}:
1\leq i \leq (s+k)(n-1)-k$ of $sF_{n+1} \cup kF_n$, the weight of
$C^{(i)}_{3}$ is $\frac{s+k}{2}(17n)+s-7k+3$. Hence $sF_{n+1} \cup
kF_n$ is $C_3$-supermagic.\qed\vspace{.5cm}

\begin{theorem}
The graph  $G\cong sL_{n+1}\cup
kL_n$ is $C_4$-supermagic, where $s, k \geq1$ and $n\geq 2$.
\end{theorem}
\noindent{\bf Proof.}\\
Let $G$ be a graph and let $v=|V(G)|$ and $e=|E(G)|$. Then
$v=2[sn+k(n-1)]$, $e=s(3n-2)+k(3n-5)$. We denote the vertex and edge
sets of $G$ as follows:

\noindent $V(G)=\{u^{j}_{i}, 1\leq i\leq n
, 1\leq j\leq s \}\cup \{ v^{j}_{i}: 1\leq i\leq n , 1\leq j\leq s
\}\cup
\{a^{t}_{i}, 1\leq i\leq n , 1\leq t\leq k \}\cup \{
b^{t}_{i}: 1\leq i\leq n , 1\leq t\leq k. \}$ \\

\noindent $E(G)=\{u^{j}_{i} v^{j}_{i}: 1\leq i\leq
n , 1\leq j\leq s \} \cup  \{u^{j}_{i}u^{j}_{i+1}:1\leq i\leq n-1 ,
1\leq j\leq s \}\cup
\{ v^{j}_{i}v^{j}_{i+1}: 1\leq i\leq n-1 ,
1\leq j\leq s \}\cup
\{a^{t}_{i} b^{t}_{i}: 1\leq i\leq n , 1\leq t\leq k \} \cup
\{a^{t}_{i}a^{t}_{i+1}:1\leq i\leq n-1 , 1\leq t\leq k \}\cup
 \{ b^{t}_{i}b^{t}_{i+1}: 1\leq i\leq n-1 ,
1\leq t\leq k \}.$\\

\noindent
 Define a total labeling $\lambda:V(G)\cup E(G)\rightarrow
\{1,2,....5n(s+k)-7k-2s\}$ as follows:\\

\noindent For $1\leq i\leq n , 1\leq j\leq s, 1\leq t\leq k.$
$$\lambda (u_i^j)=i+n(j-1)$$
$$\lambda (a_i^t)=sn+i+(n-1)(t-1)$$
$$\lambda (v_i^j)=2sn+2k(n-1)-n(j-1)-i+1$$
$$\lambda (b_i^t)=2sn+2k(n-1)+n(1-s)+t(1-n)-i$$
$$\lambda(u_i^j v_{i} ^ j)=2sn+2k(n-1)+(n-i)(s+k)+j$$
$$\lambda(a_i^t b_{i} ^ t)=2sn+2k(n-1)+(n-i-1)(s+k)+s+t$$

\noindent For $1\leq i\leq n-1 , 1\leq j\leq s, 1\leq t\leq k.$
$$\lambda(u_i^j u_{i+1} ^ j)=s(5n-2)+k(5n-7)-(s+k)(n-i-1)-j+1$$
$$\lambda(a_i^t a_{i+1} ^ t)=s(5n-2)+k(5n-7)-(s+k)(n-i-2)-s-t+1$$
$$\lambda(v_i^j v_{i+1} ^ j)=4sn+4k(n-1)-(s+k)(n-i)-j+1$$
$$\lambda(b_i^t b_{i+1} ^ t)=4sn+4k(n-1)-(s+k)(n-i-1)-s-t+1$$

\noindent
 To show that $\lambda$ is a $C_4$-supermagic labeling of $G,$
 let  $C^{(i)}_{4}, 1\leq i \leq (s+k)(n-1)-k,$ be the subcycles of
$G\cong sL_{n+1}\cup
kL_n.$ Since, the weight of every $C^{(i)}_{4}$ is $4+17n(s+k)-19k-2s,$ therefore, $G$ is a $C_4$-supermagic
\\

\section{Conclusion}
In this paper we studied the problem that if a graph $G$ has a cycle-supermagic
labeling then either  disjoint union of  isomorphic and non isomorphic copies of $G$ will have a cycle-supermagic
labeling or not? We have studied the case only for the disjoint union of $m\geq2$ isomorphic copies of fans, wheels, ladders graphs.
We also proved that disjoint union of non isomorphic copies of fans and ladders are cycle-supermagic. Moreover,
we believe that if a graph has a cycle-(super)magic labeling, then disjoint union of that graph
also has a cycle-(super)magic labeling. Therefore, we propose the following
open problem.
\begin{open}
If a graph $G$ is a cycle-(super)magic, determine whether there is a  cycle-(super)magic labeling for $mG,$ $m\geq2.$
\end{open}

\end{document}